\documentclass[a4paper,10pt]{article}
\usepackage[utf8]{inputenc}
\usepackage{amsmath,amsfonts,amssymb,amscd,amsthm,xspace, graphicx}
\usepackage[normalem]{ulem}
\usepackage[mathscr]{eucal}
\usepackage{color}
\usepackage{ulem}
\usepackage[table]{xcolor}

\usepackage[pagebackref,colorlinks,linkcolor={blue},citecolor={red}]{hyperref}
\newtheoremstyle{plain}{\topsep}{\topsep}{\slshape}{}{\bfseries}{.}{.5em}{}
\newtheorem{theorem}{Theorem}
\newtheorem{lemma}[theorem]{Lemma}
\newtheorem{corollary}[theorem]{Corollary}

\theoremstyle{remark}
\newtheorem{remark}[theorem]{Remark}
\numberwithin{equation}{section}

\def\div{\hbox{\rm div}\,}
\def\dist{\hbox{\rm dist}\,}

\def\R{{\mathbb R}}

\def\brf{\bar f}
\def\brw{\bar{\mathbf w}}

\def\me{{m}}
\def\mre{{m_0}}
\def\n{{\mathbf{n}}}

\def\we{{\mathbf w}}
\def\e{{\varepsilon}}

\def\0{{\mathbf 0}}
\def\ae{{\mathbf a}}

\def\wei{{\mathbf w}_\infty}
\def\ue{{\mathbf u}}

\def\d{\partial}

\def\XXint#1#2#3{{\setbox0=\hbox{$#1{#2#3}{\int}$ }
\vcenter{\hbox{$#2#3$ }}\kern-.6\wd0}}

\begin{document}

\title{On basic velocity estimates for the plane steady-state Navier--Stokes system and its applications }
\author{Mikhail Korobkov\footnotemark[2]  \and Xiao Ren\footnotemark[3]}
\renewcommand{\thefootnote}{\fnsymbol{footnote}}
\footnotetext[2]{
School of Mathematical Sciences, Fudan University, Shanghai 200433, P.R.China; and Sobolev Institute of Mathematics, pr-t Ac. Koptyug, 4, Novosibirsk, 630090, Russia.  email: korob@math.nsc.ru}
\footnotetext[3]{Beijing International Center for Mathematical Research, Peking University, Beijing 100871, P.R.China. email: renx@pku.edu.cn}

\maketitle

\begin{abstract}We consider some new estimates for general steady Navier-Stokes solutions in plane domains. 
According to our main result, if the domain is convex, then the difference between mean values of the velocity over two concentric circles is bounded (up to a constant factor) by the square-root of the Dirichlet integral in the annulus between the circles. The constant factor in this inequality is universal and does not depend on the ratio of the circle radii. Several applications of these formulas are discussed.
\end{abstract}

\renewcommand{\thefootnote}{\arabic{footnote}}

\section{Introduction and main results}

We study the stationary Navier--Stokes equations 
\begin{equation}  \tag{NS}\label{NS}
\left\{
\begin{aligned}
 & - \Delta \mathbf{w} + (\mathbf{w} \cdot \nabla) \mathbf{w} + \nabla p = \mathbf0 \ \ \mathrm{in} \ \Omega,\\
 & \nabla \cdot \mathbf{w} = 0 \ \ \mathrm{in} \ \Omega\\
\end{aligned}
\right.
\end{equation} 
in a plane domain (open connected set)  $\Omega\subset \R^2$,  which may be unbounded. Here $\we, p$ are the unknown velocity and pressure fields of the fluid.  
More specifically, we consider  $D$-solutions $\we$ to~(\ref{NS}) in $\Omega$, i.e., solutions with finite Dirichlet integrals 
\begin{equation} \label{unif-D-1}
\int_{\Omega} |\nabla \we|^2 < + \infty.
\end{equation}
Because of standard elliptic estimates, such solutions are $C^\infty$-smooth, and, moreover, real-analytic in $\Omega$. 

Still there are a lot of open problems concerning stationary solutions in exterior plane domains, for example, the simplest in formulation but very important from the physical point of view {\it flow around an obstacle problem},
solved for 3d case almost a century ago, but still open for 2d case:
\begin{equation}  \tag{OBS}\label{NSobs}
\left\{
\begin{aligned}
 & - \Delta \mathbf{w} + (\mathbf{w} \cdot \nabla) \mathbf{w} + \nabla p = \mathbf0 \ \ \mathrm{in} \ \Omega,\\
 & \nabla \cdot \mathbf{w} = 0 \ \ \mathrm{in} \ \Omega,\\
 & \mathbf{w}|_{\partial \Omega} = \mathbf0, \\
 & \mathbf{w} \to\mathbf{w}_\infty=\lambda \mathbf{e}_1 \ \ \text{as}\ \  |z| \to \infty.
\end{aligned}
\right.
\end{equation} 
Here $\Omega = \R^2 \setminus \bar{U}$ is an exterior plane domain, \ $U$ is the corresponding bounded open set (not necessarily connected) with smooth boundary in $\R^2$; and  $\we_\infty$ is the far field constant velocity. The parameter $\lambda > 0$ will be referred to as the Reynolds number. Here $\mathbf{e}_1=(1,0)$ is the unit vector along $x$-axis. Physically, the system \eqref{NSobs} describes the stationary motion of a viscous incompressible fluid flowing past a~rigid cylindrical body. The existence of solutions to \eqref{NSobs} with arbitrary $\lambda$ was included by Professor V.I. Yudovich in the  list of ``Eleven Great Problems in Mathematical Hydrodynamics" \cite{Yud03}. For small $\lambda$ the existence of the solutions was proved in the classical papers by R.~Finn and D.R.~Smith in 1967~\cite{FS67}, see also~\cite{KR21,KR22} for the global uniqueness of these solutions in the class of all $D$-solutions, and \cite{KR23} for the recent survey on the topic. 

The studying of $D$-solutions in unbounded plane domains is not easy at all in view of the important fact that the finiteness of the Dirichlet integral on the plane does not guarantee the boundedness of the function itself even in the mean integral sense. For instance, for  $f(z)=\bigl(\ln(2+|z|)\bigr)^{\frac13}$ we have $\int\limits_{\R^2}|\nabla f|^2<\infty$, but $f(z)\to+\infty$ as $|z|\to\infty$. This is in sharp contrast to the 3d case, where by Sobolev Embedding theorem one has
$$\biggl(\ \int\limits_{\R^3}|\nabla f|^2<\infty\biggr)\Rightarrow \biggl(\ \int\limits_{\R^3}|f-c|^6<\infty\biggr)$$ for some constant $c\in\R$. Really, the lack of the corresponding Embedding Theorem bears the main responsibility for the difficulties in the flow around an obstacle problem~(\ref{NSobs}) in 2d case and in the so called 'Stokes paradox' (see, e.g., \cite{R10}). Nevertheless, D.~Gilbarg and H.F.~Weinberger proved in~\cite{GW78} the uniform boundedness and the uniform convergence of the pressure for arbitrary $D$-solution. In particular, they established the following 
estimate.

\begin{lemma}[\cite{GW78}] \label{lem-pmean}
{\sl
Let $\we$ be a $D$-solution to the Navier--Stokes equations \eqref{NS} in the annulus $ \Omega_{r_{1},r_{2}} =  \{z \in \R^2: r_1<|z|< r_2\}$. Then for the corresponding pressure $p$ we have
\begin{equation} \label{pmean}
|\bar{p}(r_2) - \bar{p}(r_1)| \le \frac{1}{2\pi} \int\limits_{\Omega_{r_1, r_2}} |\nabla \mathbf{w}|^2
\end{equation}
 where 
 $$\bar p(r)=\frac1{2\pi r}\int\limits_{|z|=r}p(z)\,ds$$
 is the mean value of the pressure~$p$ over the circle $S_r=\{z\in\R^2:|z|=r\}$.}
\end{lemma}

The result was surprising since under general assumption 
$f \in W^{1,2}(\Omega_{r_1,r_2})$   only the simple weaker estimate 
\begin{equation} \label{bardiff}
 |\bar{f}(r_2) - \bar{f}(r_1)| \le  \sqrt{\frac{1}{2\pi}\ln\frac{r_2}{r_1}}\left(\int\limits_{\Omega_{r_1, r_2}} |\nabla{f}|^2\right)^\frac12.
\end{equation}
holds with coefficient going to $+\infty$ as $\frac{r_2}{r_1}\to\infty$. 

The inequality (\ref{bardiff}) demonstrates the well-known fact: the $D$-function in general case may have a~logarithmic growth (for example, $f(z)=\bigl(\ln(2+|z|)\bigr)^{\alpha}$ with $\alpha\in (0,\frac12)$\,). Nevertheless, the brilliant structures of Navies--Stokes system allow to obtain much better estimates. The following result was obtained just recently. 

\begin{theorem}[\cite{KR23,KR24}] \label{prop-1-basic-estimate}
{\sl Let $\we$ be the $D$-solution to the Navier--Stokes equations \eqref{NS} in the annulus ${\Omega}_{r_{1},r_{2}} $. Then we have
 \begin{equation} \label{in:estim-m1}
 |\bar\we(r_2)-\bar\we(r_1)|\le C_*\sqrt{\ln(2+\mu)}\,\biggl(\int_{\Omega_{r_1, r_2}} |\nabla \mathbf{w}|^2\biggr)^\frac12,
 \end{equation}
 where
 \begin{equation} \label{in:estim-m2}
\mu=\frac1{r_1\me},\ \quad \me:=\max\limits_{r\in[r_1,r_2]}|\bar\we(r)|,  \end{equation}
 and $C_*$ is some universal positive constant (independent of $\we, r_i$, \emph{etc.}).}
\end{theorem}
In our previous papers it was called "the first basic estimate for the velocity", its preliminary version was obtained in~\cite{GKR21}. The estimate~(\ref{in:estim-m1}) looks similar to~(\ref{bardiff}) (both contain logarithmic factor), but really it has the different nature:
$\mu(r_1,r_2)$ does {\bf not} go to $+\infty$ as $r_2\to\infty$!  \ Thus, for $r_2 \gg r_1$, \eqref{in:estim-m1} significantly improves~(\ref{bardiff}) obtained  without using the Navier-Stokes equations. Moreover, the estimate~\eqref{in:estim-m1} is precise and can not be improved in general. Indeed, for a solution 
 to~(\ref{NSobs}) in case of small $\lambda$ the opposite inequality
 \begin{equation} \label{in:estim-small-l}\int\limits_\Omega|\nabla\we|^2\le C\,\frac1{\log(2+\frac1\lambda)}\lambda^2.
  \end{equation}
 holds with $C=C(\Omega)$ (see~\cite{FS67,KR21}\,). In other words, the inequality~\eqref{in:estim-m1} has a~sharp form and the logarithmic factor there can not be omitted in general. 
 
Nevertheless, (\ref{in:estim-m1}) can be improved under some additional assumptions, for example, in the asymptotic case $r_1 \to +\infty, \frac{r_2}{r_1} \to +\infty$ \ (see \cite{GKR21}\,). The first main result of the present paper is the following 
 strengthened version of above basic estimate surprisingly having an~extremely simple form.
 
 \begin{theorem}[{\bf the improved basic velocity estimate}] \label{i-basic-estimate} {\sl
Let $\we$ be a~$D$-solution to the Navier--Stokes equations \eqref{NS} in the disk $B_{R}=\{z\in\R^2:|z|<R\}$. Then for any annulus ${\Omega}_{r_{1},r_{2}} $ with $0\le r_1<r_2\le R$ the estimate 
 \begin{equation} \label{i:estim-m1}
 |\bar\we(r_2)-\bar\we(r_1)|\le C_*\biggl(\int_{\Omega_{r_1, r_2}} |\nabla \mathbf{w}|^2\biggr)^\frac12
 \end{equation} holds, where $C_*$ is some universal positive constant (independent of $\we, r_i$, \emph{etc.}).}
\end{theorem}

In other words, logarithmic factor in (\ref{in:estim-m1}) can be dropped if we know in addition that $\we$ is a solution to Navier--Stokes system in {\bf the whole disk $B_{r_2}$}! We call the estimate~(\ref{i:estim-m1}) \ "{\it unbelievable}" because of its particular simple structure. Taking $r_2=R$, $r_1\to0$, we can obtain immediately a simple pointwise estimate:

 \begin{corollary}\label{pw-cor0} {\sl
Let $\we$ be a~$D$-solution to the Navier--Stokes equations \eqref{NS} in the disk $B_{R}=\{z\in\R^2:|z|<R\}$. Then 
 \begin{equation} \label{pw:f0}
 |\we(0)-\bar\we(R)|\le C_*\biggl(\int_{B_R} |\nabla \mathbf{w}|^2\biggr)^\frac12,
 \end{equation}
 where $C_*$ is some universal positive constant (independent of $\we,R$, \emph{etc.}).}
\end{corollary}

Really, we can improve the last assertion and obtain the similar pointwise estimates for inner subdiscs.  

 \begin{corollary}\label{pw-cor1} {\sl
Let $\we$ be a~$D$-solution to the Navier--Stokes equations \eqref{NS} in the disk $B_{R}=\{z\in\R^2:|z|<R\}$. Then for any $\delta\in(0,1)$ the estimate 
 \begin{equation} \label{pw:f1}
 \sup\limits_{z\in B_{\delta R}}|\we(z)-\bar\we(R)|\le C_\delta\biggl(\int_{B_R} |\nabla \mathbf{w}|^2\biggr)^\frac12
 \end{equation} holds, where $C_\delta$ is some universal positive constant depending on~ $\delta$ only.}
\end{corollary}

 \begin{corollary}\label{pw-cor2} {\sl
Let $\we$ be a~$D$-solution to the Navier--Stokes equations \eqref{NS} in the annulus ${\Omega}_{r_{1},r_{2}} $. Then for any $\delta\in(0,1)$ the estimates 
 \begin{equation} \label{pw:f2}
 \sup\limits_{\frac1\delta r_1\le|z|\le\delta r_2}|\we(z)-\bar\we(r_2)|\le C_{\delta}\sqrt{\ln\frac{r_2}{r_1}}\,\biggl(\int_{\Omega_{r_{1},r_{2}}} |\nabla \mathbf{w}|^2\biggr)^\frac12
 \end{equation}
  \begin{equation} \label{pw:f3}
 \sup\limits_{\frac1\delta r_1\le|z|\le\delta r_2}|\we(z)-\bar\we(r_2)|\le C_{\delta}\sqrt{\ln(2+\mu)}\,\biggl(\int_{\Omega_{r_1, r_2}} |\nabla \mathbf{w}|^2\biggr)^\frac12,
 \end{equation}
 hold, where $C_\delta$ is some universal positive constant depending on~ $\delta$ only, and
 \begin{equation} \label{in:estim-m2-cor}
\mu=\frac1{r_1\me},\ \quad \me:=\max\limits_{r\in[r_1,r_2]}|\bar\we(r)|.  \end{equation}
 }
\end{corollary}

We emphasize, that the classical Stokes estimate are far from to be enough to obtain such type inequalities (see, e.g., Remark~\ref{rem01}). 

Note, that the above estimates allow to get very simple proofs (see, e.g., \cite{KR24} for the details) of many recent advances concerning steady-state NS solutions in plane exterior domains, namely: 

(1) Every $D$-solution is uniformly bounded, and, moreover, uniformly convergent to some constant vector at infinity (see \cite{A88} under some additional assumptions, and \cite{KPR20,KPR19} for the~general case). 
\medskip

(2) Every Leray solution constructed for the problem~(\ref{NSobs}) is nontrivial~\cite{KPR21}. Moreover, it attains the required limit at infinity in case of small Reynolds numbers (=\,small~$\lambda$\,)~\cite{KR22}.

\medskip

(3) In case of small Reynolds numbers, solutions to~(\ref{NSobs}) are globally unique in the class of all $D$-solutions~\cite{KR21} .  

\medskip

Of course, (\ref{i:estim-m1}) is not only allowed to re-prove easily some previous theorems, but also to obtain some new results, for example, the following

 \begin{theorem} \label{half-plane-t}
{\sl Let $\we$ be a~$D$-solution to the Navier--Stokes equations \eqref{NS} in the domain $\Omega$ which coincides with half-plane $\R_+^2=\{z=(x,y)\in\R^2:y>0\}$ in the neighborhood of infinity:
$$\Omega=\R_+^2\setminus B_{R_0}$$
for some $R_0>0$. 
Denote $L=\partial\R^2_+$ and suppose that 
 \begin{equation} \label{eqh1}
 \lim\limits_{z\in L,\ |z|\to\infty}\we(z)=\we_0
 \end{equation} 
for some constant vector $\we_0\in\R^2$. Then the uniform convergence 
 \begin{equation} \label{eqh1-pp}
 \lim\limits_{z\in \R^2_+,\ |z|\to\infty}\we(z)=\we_0
 \end{equation}
holds.}
\end{theorem}

For comparison, in a very recent paper \cite{CLP} the convergence in similar situation\footnote{Note, that in~\cite{CLP} the authors considered stationary solutions to the Navier--Stokes system in the half-plane in presence of a~force term in the right-hand side of~\eqref{NS}.} was established only in some integral sense and only in cones separated from~$L$.

\subsection{Applications to the Leray "invading domain method"} 

Let $\Omega$ be an exterior domain in $\R^2$, i.e., 
\begin{equation}\label{Omega}
\Omega={\mathbb R}^2\setminus\bigcup_{i=1}^N\overline\Omega_i,
\end{equation}
where $\Omega_i$ are $N$  pairwise disjoint bounded Lipschitz
domains, $\overline\Omega_i\cap \overline\Omega_j=\emptyset$, \,$i\ne j$\,). 
Consider the exterior problem
\begin{equation}  \tag{GEN}\label{NSgen}
\left\{
\begin{aligned}
 & - \Delta \mathbf{w} + (\mathbf{w} \cdot \nabla) \mathbf{w} + \nabla p = 0,\\
 & \nabla \cdot \mathbf{w} = 0,\\
 & \mathbf{w}|_{\partial \Omega} = \mathbf{a}, \\
 & \mathbf{w}(z) \to\mathbf{w}_\infty=\lambda \mathbf{e}_1 \ \ \text{as}\ \  r = |z| \to +\infty,
\end{aligned}
\right.
\end{equation} 
where again $\wei$ is a constant vector from $\R^2$ and $\ae$ are boundary data (supposed to be regular enough). 
Under additional condition, that the flux of $\ae$ across each connected component of $\d \Omega$ vanishes, 
\begin{equation} \label{flux0}
\int_{\d\Omega_i} \ae \cdot \n\,ds= 0,\qquad i=1,\dots,N,
\end{equation}
where $\n$ is the unit outward (with respect to $\Omega$) normal vector, 
in the celebrated paper~\cite{L33}, J.Leray's suggested the following  elegant  approach which was called  ``{\it the~invading domains method}\,''.
Denoting by ${\mathbf w}_k$ the solution to the steady Navier--Stokes system in the  intersection of $\Omega$ with the disk $B_{R_k}$ of large  radius $R_k\ge k$ under natural boundary conditions
 \begin{equation}
\label{NSgen-k}
\left\{\begin{array}{r@{}l}
- \Delta{{\mathbf w}_k}+({\mathbf w}_k\cdot\nabla){\mathbf w}_k+ \nabla p_k  & {} ={\bf 0}\qquad\qquad\quad \hbox{\rm in } \Omega_k=\Omega\cap B_{R_k}, \\[2pt]
\nabla \cdot {{\mathbf w}_k} & {} =0\,\qquad\qquad\quad \hbox{\rm in } \Omega_k,  \\[2pt]
{{\mathbf w}_k} & {} = \ae\,\qquad\qquad\quad \hbox{\rm on } \partial\Omega,  \\[2pt]
{\mathbf w}_k & {}=\we_\infty\;\quad\mbox{for \ }|z|=R_k.
 \end{array}\right.
\end{equation}
 
 J.~Leray proved, that under~(\ref{flux0}) assumption,
the sequence ${\mathbf w}_k$ satisfies the   estimate
\begin{equation}
\label{bound}
\int_{\Omega\cap B_{R_k}}|\nabla{\mathbf w}_k|^2\le c,
\end{equation} for some positive constant $c$ independent of $k$.
 Hence, he observed that it is possible to extract a subsequence ${\mathbf w}_{k_n}$ which weakly\footnote{This convergence is uniform on every bounded set.} converges to a solution  ${\mathbf w}_L$ of problem~(\ref{NSgen})${}_{1,2,3}$ with  $\int\limits_\Omega|\nabla{\mathbf w}_L|^2<+\infty$.
This solution   was later called {\it Leray's solution\/}  (see, $e.g.$, \cite{GW74}\,). 
Recently the invading domain method was extended to the case of zero total flux: 

\begin{theorem}[\cite{KPR20}]
\label{Dbd} {\sl Let $\Omega\subset\R^2$ be an exterior domain~(\ref{Omega}) with
$C^2$-smooth boundary, and let $\wei\in\R^2$. Suppose that $\ae\in
W^{3/2,2}(\partial\Omega)$ and the equality~
\begin{equation} \label{flux-t0}
\int_{\d\Omega} \ae \cdot \n\,ds= 0
\end{equation}
holds. Then the solutions~$\we_k$ to the system~(\ref{NSgen-k}) satisfy  
the uniform estimate~(\ref{bound}). }
\end{theorem}

The proof of the last theorem is rather nontrivial and involves methods of the paper~\cite{KPR15} where the boundary value problem for~(\ref{NS}) was solved in bounded plane domains under necessary and sufficient condition of zero total flux\footnote{In fact, Theorem~\ref{Dbd} was proved in~\cite{KPR20}  only for a special case $\wei=\0$. But the proof is easily transferred to the case of an arbitrary constant vector~$\wei\in\R^2$. Such an adaptation of arguments from the case of $\wei=\0$ to an arbitrary~$\wei$ was described in details in the paper~\cite[Section~6]{KPR18} for a~more difficult case of an axisymmetric problem in an exterior 3D domain. A similar adaptation for the present plane case is carried out according to the same scheme, with many obvious simplifications.}. 

Thus, under assumption~(\ref{flux-t0}) the Leray solution $\we_L$ is still well defined and satisfies (\ref{NSgen})${}_{1,2,3}$ with  $\int_\Omega|\nabla{\mathbf w}_L|^2<+\infty$.
But the question, does the Leray solution satisfy the limiting condition at infinity
~(\ref{NSgen})${}_{4}$, is still unsolved, and it 
represents one of the significant open problems of mathematical fluid mechanics
(see \cite{KR22,GKR21} for positive answers for small $\lambda$, also \cite{KPR21} for nontriviality of Leray solutions applied to the flow around an obstacle ptoblem~(\ref{NSobs})\,). 
Thus, it seems reasonable and relevant to study the subtle properties of these solutions in order to achieve progress in understanding of their asymptotic behavior concerning the convergence~(\ref{NSgen})${}_{4}$. The next theorem of the paper represents one more step on this direction.

\begin{theorem}
\label{TBD} {\sl Under assumptions of Theorem~\ref{Dbd}, the uniform bound 
\begin{equation}
\label{bound-i}
\sup\limits_{z\in\Omega\cap B_{R_k}}|\we_k(z)|\le c\end{equation}
holds for some constant $c$ independent of~$k$.
}
\end{theorem}

Recall, in the classical paper~\cite{GW74} D.~Gilbarg and H.F.~Weinberger proved that 
for any fixed $\delta\in (0,1)$ the uniform estimates 
\begin{equation}
\label{FD3-nn}\sup\limits_{\{z\in\Omega: |z|<\delta R_k\}}|\we_k(z)|\le C=C(\delta),
\end{equation}
\begin{equation}
\label{FD4-nn}\sup\limits_{\{z\in\Omega: |z|<\delta R_k\}}|p_k(z)|\le C=C(\delta)
\end{equation}
hold. 
Of course, (\ref{bound-i}) is much stronger than~(\ref{FD3-nn}). Moreover, (\ref{bound-i}) allows us to obtain some new interesting estimates for Leray solutions. 
Here we formulate the simplest variant only. 

Consider the case when $\Omega$ and $\mathbf{a}=(a_1,a_2)$ are symmetric with respect to the $x$-axis. Precisely, it is assumed that
\begin{equation} \label{eq-Omegasym}
(x, y) \in \Omega \iff (x, -y) \in \Omega
\end{equation}
and
\begin{equation} \label{eq-asym}
a_1(x, y) = a_1(x, -y), \quad a_2(x, y) = -a_2(x, -y).
\end{equation}

\begin{theorem}\label{STRbd}
{\sl Under assumptions of Theorem~\ref{Dbd}, suppose, in addition, that $\we_\infty=\mathbf{0}$, $\partial\Omega\in C^3$,  and that $\Omega$ and the boundary data $\mathbf{a}\in W^{5/2,2}(\Omega)$ satisty the conditions \eqref{eq-Omegasym}--\eqref{eq-asym}. Take the sequence of symmetric solutions~$\we_k$ to~(\ref{NSgen}). 
Then the estimate
\begin{equation}
\label{gg1} \int\limits_{\Omega\cap B_{R_k}}\psi_k^2|\nabla\omega_k|^2\le C
\end{equation}
holds
for some constant~$C$ which is independent of~$k$, where $\psi_k$ are the corresponding stream functions
$$\nabla\psi_k=\we_k^\bot=(-w_{k2},w_{k1}),
\qquad\psi_k|_{S_{R_k}}\equiv0,$$
and $\omega_k=\Delta\psi_k=\partial_yw_{k1}-\partial_xw_{k2}$ is the corresponding vorticity.}
\end{theorem}
Here, a solution $\we=(w_1,w_2)$ is said to be symmetric if
\begin{equation} \label{eq-wsym}
w_1(x, y) = w_1(x, -y), \quad w_2(x, y) = -w_2(x, -y).
\end{equation}

We emphasize, that the previous estimates (\ref{FD3-nn})--(\ref{FD4-nn}) do  {\bf not} allow us to prove boundedness of the integral~(\ref{gg1}) even in subdomains $\Omega\cap B_{\delta R_k}$ with fixed $\delta\in(0,1)$.

The results of the type of Theorems~\ref{TBD}--\ref{STRbd}  can presumably shed new light on the nature of Leray's solutions and contribute to their further study.

\medskip

Let us say a few words concerning proofs. The proof of the previous basic estimate in Theorem~\ref{prop-1-basic-estimate} was based on some delicate analysis of Bernoulli pressure level sets, see, e.g.,  \cite{KPR21,GKR21} for the description of the arguments on a~heuristic level.  The main novelty in the proof of the improved version in Theorem~\ref{i-basic-estimate} (without logarithmic factor) is the using of the maximum principle for the vorticity~$\omega$ in the whole disk~$B_R$ (see below Corollary~\ref{lem:Stokes3}). Roughly speaking, if $r_1>\frac1m$, then we can apply the initial Theorem~\ref{prop-1-basic-estimate} with trivially bounded logarithmic factor, but if $r_1<\frac1m$, we apply another inequality using the vorticity uniform estimates in the half-disk (see~(\ref{est-sf6}) and Corollary~\ref{lem:Stokes4}\,). The unexpected and exciting fact is, that these two quite different approaches are in mutual harmony with each other, and finally allow us to get the required estimate~(\ref{i:estim-m1}) of an~amazingly simple form.
 
The paper is organized as follows: we give some preliminary technical in the next section, and the proof of the above results are contained in section~\ref{proofs-sec}.

\section{Preliminaries} \label{sec-fs}
 
 \subsection{Some properties of $D$-functions} \label{sec:2.0}
 For a function $f(z)$ defined on a domain $\Omega\subset\R^2$ we denote by 
 $$\int\limits_\Omega f(z)$$  the integral with respect to the usual two-dimensional Lebesgue measure, and by 
  $$\int\limits_\Gamma f(z)\,ds$$
  the integral with respect to the curve length (=$1$-dimensional Hausdorff measure). 
 
We use standard notations for Sobolev spaces $W^{1,2}(\Omega)$. 
Below for $\rho_2>\rho_1>0$ we denote by $\Omega_{\rho_1, \rho_2}$ the corresponding annulus-type domain:
$$\Omega_{\rho_1, \rho_2} = \{z \in \mathbb{R}^2: \rho_1 <|z| < \rho_2\}.$$
Further, for a function $f\in W^{1,2}({\Omega}_{\rho_1, \rho_2})$ denote by $D_f(\rho_1, \rho_2)$ the corresponding 
Dirichlet integral:
$$D_f(\rho_1, \rho_2):=\int\limits_{\Omega_{\rho_1, \rho_2}}|\nabla f|^2.$$
Using this notation, recall three elementary facts concerning $D$-functions. 

\begin{lemma} \label{lem-new4}
{\sl For any $f\in W^{1,2}({\Omega}_{\rho_1, \rho_2})$ we have 
\begin{equation} \label{eq-new2.9}
|\bar{f}(\rho_2) - \bar{f}(\rho_1)| \le \frac{1}{\sqrt{2\pi}}  \left(\ln \frac{\rho_2}{\rho_1}\right)^\frac12  \cdot\bigl(D_f(\rho_1,\rho_2)\bigr)^\frac12,
\end{equation}
where $\brf$ means the mean value of~$f$ over the circle $S_r$:
$$\brf(r):=\frac1{2\pi r}\int\limits_{|z|=r}f(z)\,ds.$$}
\end{lemma}

 \begin{lemma}
\label{lt2} {\sl Fix a number $\beta\in(0,1)$. Let $R>0$ and $f\in W^{1,2}(\Omega_{\beta R, R} )$. Then there exists  a number $r\in[\beta R, R]$ such that
 the
estimate \begin{equation}
\label{est-2}
\sup\limits_{|z|=r}|f(z)-\brf(r)|\le c_\beta \sqrt{D_f(\beta R,R)} 
\end{equation} holds, where~$c_\beta=\sqrt{\frac{2\pi}{1-\beta}}$.   }
 \end{lemma}
 
 Summarize the results of these lemmas, we receive
\begin{lemma}\label{lt3}
{\sl Under conditions of Lemma~\ref{lt2}, there exists $r\in [\beta R,R]$ such that
\begin{equation}
\label{est-3}
\sup\limits_{|z|=r}|f(z)-\brf(R)|\le \tilde c_\beta\sqrt{D_f(\beta R,R)}, 
\end{equation} holds, where~the positive constant $\tilde c_\beta$ depends on $\beta$ only.  
}
\end{lemma}

The circles $S_r$ from lemmas~\ref{lt2}--\ref{lt3} satisfying the formulated uniform pointwise estimates (\ref{est-2})--(\ref{est-3}) are usually called ''good'' circles; this simple approach is often used in the studying of $D$-solutions. 

\subsection{The Stokes estimates} \label{sec:2.3}
We recall the following simplest local regularity estimate for the linear Stokes system. For the proof, see, for instance, \cite[Theorem IV.4.1 and Remark IV.4.1]{G11}. 

\begin{lemma} \label{lem:Stokes}{\sl 
Let $\mathbf{w}_S$ be a local solution in the annulus $\Omega_{1,4}=\{z\in\R^2:1<|z|<4\}$ to the Stokes system
\begin{equation} \label{Stokes}
\left\{
\begin{aligned}
 & \Delta \mathbf{w}_S - \nabla p_S = \mathbf{f}_S, \\
 & \nabla \cdot \mathbf{w}_S = 0. \\
\end{aligned}
\right.
\end{equation}
Then there holds the following regularity estimates for  $1<s<\infty$:
\begin{align} \label{StokesEstimate}
 \|\nabla^{2} \mathbf{w}_S\|_{L^s(\Omega_{2,3})} \le C(s) \,\bigl(\|\mathbf{w}_S\|_{W^{1,s}(\Omega_{1,4})} +  \|\mathbf{\mathbf{f}}_S\|_{L^s(\Omega_{1,4})}\bigr).
\end{align}}
\end{lemma}

Applying the above lemma to the partial case $s=4/3$ and using the standard scaling and the simplest Sobolev Imbedding estimates, we obtain the following useful

\begin{lemma} \label{lem:Stokes0}{\sl 
Let $\mathbf{w}$ be a $D$-solution to the Navier--Stokes system~(\ref{NS}) in an annulus $\Omega_{\frac14R,R}\subset\R^2$. Denote 
\begin{equation}
\label{est-s1}
\mre:=\brw(R),\qquad D=\int\limits_{\Omega_{\frac14R,R}}|\nabla\we|^2,
\end{equation}
and suppose in addition that 
\begin{equation}
\label{est-s2}
\sqrt{D}\le \mre.
\end{equation}
Then the estimates 
\begin{equation} \label{est-sf1}
 \|\nabla^{2} \mathbf{w}\|_{L^{4/3}(\Omega_{\frac12R,\frac34 R})}\le C  R^\frac12\sqrt{D}\biggl(\mre+\frac1{R}\biggr),
\end{equation}
\begin{equation} \label{est-sf2}
\|\nabla^{2} \mathbf{w}\|_{L^1(\Omega_{\frac12R,\frac34 R})}\le C  \sqrt{D}\bigl(1+\mre R\bigr),
\end{equation}
hold, where $C $ is a some universal constant (independent of $\we, R$, etc.). }
\end{lemma}

Of course, (\ref{est-sf2}) is an~elementary consequence of~(\ref{est-sf1}) and the Holder inequality. 

Under conditions of Lemma~\ref{lem:Stokes0}, using~(\ref{est-s1})${}_2$ and (\ref{est-sf2}) it is easy fo find a ''good'' circle $S_\rho$ with $\rho\in\bigl(\frac12R,\frac34 R\bigr)$, such that 
\begin{equation} \label{est-sf4}
\int\limits_{S_\rho}\omega^2\,ds\le \frac{8}{R}D,\qquad \int\limits_{S_\rho}|\nabla\omega|\,ds\le C \frac1R \sqrt{D}\bigl(1+\mre R\bigr),
\end{equation}
where $\omega(z)=\partial_yw_{1}-\partial_xw_{2}$ is the corresponding vorticity.  Therefore, 
$$\max\limits_{z\in S_\rho}|\omega(z)|\le C \frac1R \sqrt{D}\bigl(1+\mre R\bigr).$$ Since $\omega$ satisfies the maximum principle~(see, e.g., \cite{GW78}), we obtain the following 
useful

\begin{corollary} \label{lem:Stokes3}{\sl 
Under assumptions of Lemma~\ref{lem:Stokes0}, if an addition $\we$ is a $D$-solution to the Navier--Stokes system in the whole disk $B_R$, then the estimate 
\begin{equation} \label{est-sf5}
\max\limits_{z\in B_{\frac12 R}}|\omega(z)|\le C  \frac1R \sqrt{D}\bigl(1+\mre R\bigr)
\end{equation}
is valid, where as above $ D=\int\limits_{\Omega_{\frac14R,R}}|\nabla\we|^2$ and  $C $ is a some universal constant (independent of $\we, R$, etc.). }
\end{corollary}

For a point $z=(x,y)\in\R^2$ denote $z^\bot=(x,y)^\bot=(-y,x)$.  By the classical elegant formula from complex analysis, for any divergence-free $D$-function $\we(z)$ in annulus 
$\Omega_{r_1,r_2}$ we have
\begin{equation} \label{est-sf6}
\brw(\rho_2)-\brw(\rho_1)=\frac1{2\pi}\int\limits_{\Omega_{\rho_1,\rho_2}}\frac{\omega(z)\,z^\bot}{|z|^2}
\end{equation}
(see, e.g., \cite[p.~388]{GW78} or \cite[p.206]{Sa99}). In particular, using standard formula for two dimensional integral in polar coordinate system, we have
\begin{equation} \label{est-sf7}
\bigl|\brw(\rho_2)-\brw(\rho_1)\bigr|\le(\rho_2-\rho_1)\cdot\max\limits_{z\in \Omega_{\rho_1,\rho_2}}|\omega(z)|.
\end{equation}

Combining it with Corollary~\ref{lem:Stokes3}, we get 

\begin{corollary} \label{lem:Stokes4}{\sl 
Under assumptions of Lemma~\ref{lem:Stokes0}, if an addition $\we$ is a $D$-solution to the Navier--Stokes system in the whole disk $B_R$, and  
\begin{equation} \label{est-sf8}
0\le \rho_1< \rho_2\le R\qquad\mbox{ with }\qquad
|\rho_2-\rho_1|\le \frac1\mre,\end{equation}
then 
\begin{equation} 
\label{est-sf10}
\bigl|\brw(\rho_2)-\brw(\rho_1)\bigr|\le C \sqrt{D},
\end{equation}
where as above $ D=\int\limits_{\Omega_{\frac14R,R}}|\nabla\we|^2$ and  $C $ is a some universal constant (independent of $\we, R$, etc.). }
\end{corollary}

Indeed, if $\rho_2\le\frac12 R$, then from~(\ref{est-sf5}), (\ref{est-sf7}), (\ref{est-sf8}) we obtain directly
$$\bigl|\brw(\rho_2)-\brw(\rho_1)\bigr|\le C \sqrt{D}\min\biggl\{\bigl(1+\mre R), \frac{(1+\mre R)}{\mre R}\biggr\}\le  C' \sqrt{D}.$$
 The opposite case when $\rho_2\in\bigl[\frac12R,R\bigr]$ can be reduced to the previous one using the trivial estimate $\max\limits_{r\in[\frac12R,R]}\bigl|\brw(R)-\brw(r)\bigr|\le \sqrt{D}$ obtained from (\ref{bardiff}). 

\begin{remark}\label{rem01}
It is well-known that the $L^1$-norm of the second derivative is enough to estimate the $C$-norm of the function itself (see, e.g., \cite{BKK}), in particular,  (\ref{est-sf2}) implies
\begin{equation} \label{est-sfr}
\|\mathbf{w}-\brw(R)\|_{L^\infty(\Omega_{\frac12R,\frac34 R})}\le C' \sqrt{D}\bigl(1+\mre R\bigr).
\end{equation}
The right-hand side goes to infinity as $R\to\infty$, so this estimate is much worse than~(\ref{pw:f2}) based on more subtle tools, see below. 
\end{remark}

\section{Proofs}\label{proofs-sec}

Below, unless otherwise specified, we use $C$ to denote some universal constants (independent of $\we$, $\Omega$, $r_i$, etc.)\footnote{The obvious exception concerns the proof of Theorem~\ref{TBD} and Corollary~\ref{STRbd}. The constants $C$ used theret, of course, depend on $\Omega$ and $\ae$, but they do not depend on~$k$, this is the crucial issue there.}. The exact values of $C$ may change from line to line.

{\bf\large Proof of Theorem~}\ref{i-basic-estimate}. Let the assumptions of the theorem be fulfilled, i.e, let $\we$ be a $D$-solution to the Navier--Stokes system in the disk $B_R$. For a given annulus subdomain  $\Omega_{r_1,r_2}=\{z\in\R^2:r_1<|z|<r_2\}\subset B_R$ 
denote 
$$D(r_1,r_2)=\int\limits_{\Omega_{r_1,r_2}}|\nabla\we|^2.$$ 
So we have to prove that 
 \begin{equation} \label{f-add-t1}
\bigl|\brw(r_2)-\brw(r_1)\bigr|\le C_*\sqrt{D(r_1,r_2)} 
\end{equation}
 Let us make several simplifications. First of all, for 
it is sufficient to consider the case
 \begin{equation} \label{f-add1}
|\bar\we(r)|\ge\frac14 m\qquad \forall r\in[r_1,r_2], 
\end{equation}
where, recall, 
$$\me=\max\limits_{r\in[r_1,r_2]}|\bar\we(r)|.$$
Indeed, suppose that (\ref{f-add1}) fails. Then take $r_0\in[r_1,r_2]$ with 
$$|\bar\we(r_0)|=\me,$$ and take also  $r'_0\in [r_1,r_2]$ with 
$$|\bar\we(r'_0)|<\frac14\me.$$
Assume, for definiteness, that $r_0<r'_0$ (the case of opposite inequality can be considered analogously). 
Then we put
$$r'_1=r_0,\qquad r'_2=\min\limits_{r\in[r_0,r'_0]}\{r:|\bar\we(r)|\le\frac14\me\}.$$
 Now we have to consider the interval $[r'_1,r'_2]$ instead of $[r_1,r_2]$. 
 By construction, 
 $$\!\!\!D(r'_1,r'_2)<D(r_1,r_2),\quad\bigl|\bar\we(r_2)-\bar\we(r_1)\bigr|\le 2\me=\frac83\cdot\frac34\me=\frac83\,\bigl|\bar\we(r'_2)-\bar\we(r'_1)\bigr|.$$
 So if we prove the required estimate (\ref{f-add-t1}) for $r'_1,r'_2$, then it implies the required estimate for initial pair $r_1,r_2$ (with slightly modified constant $C_*'=\frac83 C_*$). So below we assume that~(\ref{f-add1}) is fulfilled.

Of course, we can assume also that 
 \begin{equation} \label{f-add00}
\begin{array}{rcl}
r_2=R,\\[4pt]
r_1<\frac14R,
\\[4pt]
\sqrt{D(r_1,r_2)}<\frac14 m.\end{array}
\end{equation}
Indeed, if the first equality fails, we can simply put $R'=r_2$, and consider the disk $B_{R'}$ instead of initial $B_R$. Further, if $r_1\ge\frac14 R$, then the desired inequality~(\ref{f-add-t1}) follows immediately from the trivial estimate~(\ref{bardiff}). Furthermore, 
if $\sqrt{D(r_1,r_2)}\ge\frac14 m$, then by construction $\bigl|\bar\we(r_2)-\bar\we(r_1)\bigr|\le 2\me\le 8 \sqrt{D(r_1,r_2)}$, and there are nothing to prove. So below we assume, that all assumptions in~(\ref{f-add00}) are fulfilled as well. In particular, from (\ref{f-add00})${}_{1,3}$
 and (\ref{f-add1}) we obtain   
  \begin{equation} \label{f-add01}
\sqrt{D(r_1,r_2)}<|\bar\we(R)|.
\end{equation}
Now we have to consider several possible cases.

\medskip
{\sc\large Case I: $r_1\me\ge\frac14$.} Then the required inequality~(\ref{f-add-t1}) follows immediately from the first basic estimate~(\ref{in:estim-m1}).

\medskip
{\sc\large Case II: $r_1\me<\frac14$.} Following the notation of subsection~\ref{sec:2.3}, denote
$$\mre=|\we(R)|,\qquad D=\int\limits_{\Omega_{\frac14R,R}}|\nabla\we|^2.$$
By construction,
\begin{equation} \label{ff-1}
\frac14\me\le\mre\le\me, \qquad D<D(r_1,r_2)\le\mre
\end{equation}
(see (\ref{f-add01}), (\ref{f-add00})${}_{1-2}$\,). 
Put $r_*=\frac1{4\me}$. By above 'Case II' assumption and by (\ref{ff-1})${}_{1}$, 
\begin{equation} \label{ff-2}
r_1<r_*<\frac1{\mre}.
\end{equation}
Now we have to consider two different subcases. 

\medskip
{\sc\large Case IIa: $r_*\ge R$.} This implies 
$$R\le\frac1\mre.$$ Then the required inequality~(\ref{f-add-t1}) follows immediately from Corollary~\ref{lem:Stokes4}. 

\medskip
{\sc\large Case IIb: $r_*<R$.} Then from the first basic estimate~(\ref{in:estim-m1}) applied to the annulus $\Omega_{r_*,R}$ we obtain immediately that 
$$\bigl|\brw(R)-\brw(r_*)\bigr|\le C \sqrt{D(r_*,R)}<C \sqrt{D(r_1,R)}.$$ 
Further, by Corollary~\ref{lem:Stokes4} applied to the radii  $ \rho_1=r_1,\ \rho_2=r_*$, we obtain 
$$\bigl|\brw(r_*)-\brw(r_1)\bigr|\le C \sqrt{D}<C_*\sqrt{D(r_1,R)}.$$
The last two display formulas imply the required inequality~(\ref{f-add-t1}). Theorem~\ref{i-basic-estimate} is proved completely. \hfill $\qed$

\

{\bf\large Proof of Theorem~}\ref{TBD}. Let assumptions of Theorems~\ref{Dbd}--\ref{TBD} be fulfilled. So we have a sequence of solutions 
\begin{equation}
\label{NSDI}
 \left\{\begin{array}{rcl}
-\nu \Delta{\bf w}_{k}+\big({\bf w}_{k}\cdot \nabla\big){\bf
w}_{k} +\nabla p_{k} &
= & { 0}\qquad \hbox{\rm in } \Omega_k,\\[4pt]
\div\,{\bf w}_{k} &  = & 0  \qquad \hbox{\rm in } \Omega_k,
\\[4pt]
 {\bf w}_{k} &  = & {\bf a}
 \qquad \hbox{\rm on }\partial\Omega,\\[4pt]
  {\bf w}_{k} &  = & \wei
\,\quad \hbox{\rm on }\partial B_k,\end{array}\right.
\end{equation}
in domains bounded domain
$\Omega_k=\Omega\cap B_k$, \ $B_k=\{z:|z|<R_k\}$, $R_k\to+\infty$.

Below in this proof we use $C$ to denote constants that are independent of $k$. The exact values of $C$ may change from line to line.  
By Theorem~\ref{Dbd},
\begin{equation}
\label{FDIk} \int\limits_{\Omega_k}|\nabla{\bf w}_k|^2\le C.
\end{equation}
Also, by D.~Gilbarg and H.F.~Weinberger results~\cite{GW74},  we have 
\begin{equation}
\label{FD3}\sup\limits_{\{z\in\Omega: |z|<\frac34 R_k\}}|\we_k(z)|\le C.\end{equation}

The proof is divided into the three steps. 

{\sc Step~1}. Since
\begin{equation}
\label{P1}\int\limits_{B_{k}\setminus\frac12B_k}|\nabla\we_k|^2
=\int\limits_0^{2\pi}\biggl(\ \ \int\limits_{\frac12R_k}^{R_k}r|\nabla\we_k(r,\theta)|^2\,dr\biggr)\,d\theta\le C,
\end{equation}
we have 
\begin{equation}
\label{P2}\exists\theta_k\in(0,2\pi):\qquad\int\limits_{\frac12R_k}^{R_k}|\nabla\we_k(r,\theta_k)|^2\,dr
\le\frac{C}{R_k},
\end{equation}
therefore,
\begin{equation}
\label{P3}\int\limits_{\frac12R_k}^{R_k}|\nabla\we_k(r,\theta_k)|\,dr
\le C.
\end{equation}
Hence because of the boundary conditions on~$S_{R_k}=\partial B_k$, we have
\begin{equation}
\label{P4}\sup\limits_{\frac12R_k\le r\le R_k}|\we_k(r,\theta_k)|
\le C.
\end{equation}

{\sc Step~2}. Put $\ue_k(z):=\we_k(z)-\wei$, in particular, 
\begin{equation}
\label{FD4}\ue_k(z)\equiv \0\quad\mbox{ on }S_{R_k}=\partial B_k.
\end{equation}
For $h>0$ denote 
$$\Omega_{kh}=\{z\in\R^2:R_k-h<|z|<R_k\},$$
$$S_{kh}=\{z\in\R^2:|z|=R_k-h\}$$
(below we always assume that $0<h<\frac12R_k$).

Using the classical Hardy inequality and zero boundary conditions on $S_{R_k}$, we obtain 
\begin{equation}
\label{SFD1}\int\limits_{\Omega_{kh}}|\ue_k|^2\le C\,h^2.
\end{equation}
Consequently, by Holder inequality, 
\begin{equation}
\label{H1}\int\limits_{\Omega_{kh}}|\ue_k|\cdot|\nabla\ue_k|=\int\limits_0^h\biggl(\int\limits_{S_{k\tau}}|\ue_k|\cdot|\nabla\ue_k|\,ds\biggr)\,d\tau\le C\,h. 
\end{equation}
Thus, by elementary mean-value theorem, 
\begin{equation}
\label{H2}\forall h\in\bigl(0,\frac12R_k\bigr)\ \exists \tau\in\bigl(\frac23h,h\bigr)\quad\mbox{ such that \ }\ \ \int\limits_{S_{k\tau}}|\ue_k|\cdot|\nabla\ue_k|\,ds\le C. 
\end{equation}
From the last assertion and from~(\ref{P4}) we obtain immediately the following important property:
\begin{equation}
\label{H3}\forall h\in\bigl(0,\frac12R_k\bigr)\ \exists \tau\in\bigl(\frac23h,h\bigr)\quad\mbox{ such that \ }\ \ \sup\limits_{z\in S_{k\tau}}
|\ue_k(z)|\le C.
\end{equation}

{\sc Step~3}. Take and fix arbitrary $z_0\in B_k\setminus \frac23 B_k$ and denote $r_0=|z_0|$, \ $h_0=R_k-r_0$. By construction,  
\begin{equation}
\label{Hz} \frac23 R_k<r_0<R_k,\qquad 0<h_0<\frac13 R_k.
\end{equation}
Applying (\ref{H3}) to $h=h_0$, we obtain the existence of $\tau_0\in\bigl(\frac23h_0,h_0\bigr)$ such that 
\begin{equation}
\label{H3-1}\sup\limits_{z\in S_{k\tau_0}}
|\ue_k(z)|\le C.
\end{equation}
Denote by $S_{r}(z_0)$ the circle of radius~$r$ centered at $z_0$, and by the $B_r(z_0)$ the corresponding disk $B_r(z_0):=\{z\in\R^2:|z-z_0|<r\}$. 
By construction, 
\begin{equation}
\label{H3-2}B_{h_0}(z_0)\subset\Omega_k,
\end{equation}
\begin{equation}
\label{H3-3}S_r(z_0)\cap S_{k\tau_0}\ne\emptyset\qquad\forall r\in\bigl(\frac13h_0,h_0\bigr).
\end{equation}
Applying (\ref{FDIk}) and Lemma~\ref{lt2}, we found ''good'' circle $S_{r_*}(z_0)$ with $r_*\in\bigl(\frac13h_0,h_0\bigr)$ such that 
\begin{equation}
\label{H3-4}
\sup\limits_{z_1,z_2\in S_{r_*}(z_0)}|\ue_k(z_1)-\ue_k(z_2)|\le C. 
\end{equation}
Then from (\ref{H3-1}) and (\ref{H3-3}) we finally get
\begin{equation}
\label{H3-5}
\sup\limits_{z\in S_{r_*}(z_0)}|\ue_k(z)|\le C. 
\end{equation}
Now "unbelievable" estimate~(\ref{pw:f0}) applying for the disk $B_{r_*}(z_0)$ immediately implies 
 \begin{equation}
\label{H3-6}
|\ue_k(z_0)|\le C.
\end{equation}
Because of arbitrariness of $z_0\in B_k\setminus \frac23 B_k$, we proved really that 
\begin{equation}
\label{H3-7-00}
\sup\limits_{z\in B_k\setminus \frac23 B_k}|\we_k(z)|\le C. 
\end{equation}
Together with~(\ref{FD3}) this implies the required uniform estimate
\begin{equation}
\label{H3-7}
\sup\limits_{z\in \Omega_k}|\we_k(z)|\le C. 
\end{equation}
\qed

\medskip
{\bf \large Proof of Theorem~}\ref{STRbd}. Let assumptions of Theorems~\ref{Dbd}--\ref{STRbd} be fulfilled.  In particular, by (\ref{bound}),
\begin{equation}
\label{bound-h8}
\int_{\Omega_k}\omega_k^2\le 2\int_{\Omega_k}|\nabla\we_k|^2\le C,
\end{equation}Then 
$$\Delta\omega_k=\we_k\cdot\nabla\omega_k,$$
\begin{equation}
\label{bound-h9}
\begin{aligned}
 & \div(\psi_k^2\omega_k\nabla\omega_k)=\psi_k^2|\nabla\omega_k|^2+2\psi_k\omega_k\nabla\psi_k\nabla\omega_k+\psi_k^2\omega_k(\we_k\cdot\nabla\omega_k)\\
 & =
\psi_k^2|\nabla\omega_k|^2+2\psi_k\nabla\psi_k\omega_k\nabla\omega_k+\frac12\div(\psi_k^2\omega^2_k\we_k).\\
\end{aligned}
\end{equation}
Now, integrating the identity (\ref{bound-h9}) over the domain~$\Omega_k$, using zero boundary conditions $\psi_k\equiv0$ on the big circle $S_{R_k}$, we obtain 
\begin{equation}
\label{bound-h10}
\int\limits_{\partial\Omega}\psi_k^2\omega_k\partial_n\omega_k\,dS-\frac12\int\limits_{\partial\Omega}\psi_k^2\omega^2_k\bigl(\ae\cdot\n\bigr)\,dS
=\int\limits_{\Omega_k}\psi_k^2|\nabla\omega_k|^2+2\int\limits_{\Omega_k}\psi_k\omega_k\nabla\psi_k\nabla\omega_k.
\end{equation}
Because of standard regularity properties of steady Stokes and Navier--Stokes system in smooth bounded domains (see, e.g., \cite{G11,KPR24}), from~(\ref{bound-h8}) and the assumption~$\ae\in W^{5/2,2}(\Omega)$ we obtain that 
the left hand side of the last identity is uniformly bounded (with respect to~$k$\,):
\begin{equation}
\label{bound-h11}
\int\limits_{\partial\Omega}\psi_k^2\omega_k\partial_n\omega_k\,dS-\frac12\int\limits_{\partial\Omega}\psi_k^2\omega^2_k\bigl(\ae\cdot\n\bigr)\,dS
\le C.
\end{equation}
On the other hand, from the elementary Holder inequality 
$$|\psi_k\omega_k\nabla\psi_k\nabla\omega_k|=|\psi_k\omega_k\we_k^\bot\cdot\nabla\omega_k|\overset{(\ref{H3-7})}\le C |\psi_k\omega_k\nabla\omega_k|\le  \e \psi_k^2|\nabla\omega_k|^2+\frac{C}\e\omega_k^2,$$
and from (\ref{bound-h10})--(\ref{bound-h11}) we obtain finally the required estimate
\begin{equation}
\label{bound-h12}
\int\limits_{\Omega_k}\psi_k^2|\nabla\omega_k|^2\le C.
\end{equation}
 \hfill $\qed$
 
 \medskip
 
 Note, that the standard D.~Gilbarg and H.F.~Weinberger estimates (\ref{FD3-nn})--(\ref{FD4-nn}) are {\bf not} enough to establish the estimate of kind~(\ref{bound-h12}), even in the subdomain $\Omega\cap B_{\frac12R_k}$. Here we have used the more precise estimate of Theorem~\ref{TBD} in the whole disk in essential way. 

\medskip
{\bf \large Proof of Theorem~}\ref{half-plane-t}. Let assumptions of Theorems~\ref{half-plane-t} be fulfilled. Without loss of generality we assume that  
$$\Omega=\R^2_+\setminus B_1.$$
Denote $S_R^+=S_R\cap\Omega$,  and 
$$\e^2_k=\int\limits_{\{z\in\Omega:2^k<|z|<2^{k+1}\}}|\nabla\we|^2.$$ Since $\we$ is a $D$-solution, we have 
$\e_k\to 0$ as $k\to\infty$. Take an increasing sequence $R_k\to\infty$ such that $R_k\in[2^k,2^{k+1}]$  and
\begin{equation}
\label{hpf1}
\int\limits_{S^+_{R_k}}|\nabla\we|^2\,ds\le\frac4{R_k}\e^2_k,
\end{equation}
the existence of such sequence follows from the mean-value theorem. Then by Holder inequality,
\begin{equation}
\label{hpf2}
\int\limits_{S^+_{R_k}}|\nabla\we|\,ds\le4\e_k.
\end{equation}
From the last estimate and from the assumption~(\ref{eqh1}) we have 
\begin{equation}
\label{hpf3}
\sup\limits_{z\in S^+_{R_k}}|\we(z)-\we_0|\le\delta_k,
\end{equation}
where $\delta_k:=4\e_k+\tilde\delta_k$ with $\tilde\delta_k:=\sup\limits_{z\in L, |z|\ge R_k}|\we(z)-\we_0|\to0$ as $k\to\infty$. 
Denote
$$\Omega_k=\Omega\cap B_{R_{k+1}}\setminus B_{R_k}.$$ 
By construction,
$$\partial\Omega_k=S^+_{R_k}\cup S^+_{R_{k+1}}\cup \bigl(L\cap\bar\Omega_k\bigr),$$
and 
\begin{equation}
\label{hpf3-}
\sup\limits_{z\in \partial\Omega_k}|\we(z)-\we_0|\le\delta_k+\delta_{k+1}\to0\qquad\mbox{ as }k\to\infty. 
\end{equation}

Now take and fix arbitrary a point $z_0\in \Omega_k$ and denote $r_0=\dist(z_0,\partial\Omega_k):=\inf\limits_{z\in \partial\Omega_k}|z-z_0|$. 
By construction,
$$B(z_0,r_0)\subset\Omega_k.$$
\medskip 
Further we have to consider separately three
possible cases.

($i$) \ $r_0=\sup\limits_{z\in S_{R_k}}|z-z_0|$, i.e., the point $z_0$ is near the half-circle $S^+_{R_k}$ (or, more precisely, it is closer to $S^+_{R_k}$ than to other parts of~$\partial\Omega_k$);

($ii$) \ $r_0=\sup\limits_{z\in S_{R_{k+1}}}|z-z_0|$, i.e., the point $z_0$ is near the half-circle $S^+_{R_{k+1}}$;

($iii$) \ $r_0=\sup\limits_{z\in L}|z-z_0|$, i.e., the point $z_0$ is near the horizontal axes $L\cap\bar\Omega_k$.

\medskip

{\sc Case ($i$)}: \  $r_0=\sup\limits_{z\in S_{R_k}}|z-z_0|$. For a point $z\in \Omega_k$ with polar coordinates $z\sim (r,\theta)$ define 
$$\ue_k(r,\theta)=\we(r,\theta)-\we(R_k,\theta)$$
Then by one-dimensional Hardy inequality we have
$$\int\limits_{\Omega_{0k}}|\ue_k|^2\le C\, r_0^2\e_k^2,$$
where $C$ is some universal constant, and $\Omega_{0k}:=\{z\in\Omega_k:R_k<|z|<R_k+2r_0\}$. Then there exists a ''good'' radius $r_*\in \bigl(\frac12 r_0,r_0\bigr)$ such that 
\begin{equation}
\label{hpf4}
\int\limits_{S_{r_*}(z_0)}|\ue_k|^2\,ds\le C\, r_0\e_k^2,
\end{equation}
\begin{equation}
\label{hpf5}\int\limits_{S_{r_*}(z_0)}|\nabla\we|^2\,ds\le \frac{C}{r_0}\e_k^2,
\end{equation}
consequently,
\begin{equation}
\label{hpf6}
\int\limits_{S_{r_*}(z_0)}|\nabla\we|\,ds\le C\e_k.
\end{equation}
From (\ref{hpf4}) we conclude, that there exists a point $z_*\in S_{r_*}(z_0)$ such that 
$$|\ue_k(z_*)|\le C\e_k.$$
Then by construction and from (\ref{hpf3}) we get 
$$|\we(z_*)-\we_0|\le\delta_k+C\e_k.$$
From the last estimate and from (\ref{hpf6}) we conclude
\begin{equation}
\label{hpf7}
\sup\limits_{z\in S_{r_*}(z_0)}|\we(z)-\we_0|\le\delta_k+C\e_k. 
\end{equation}
Then from the "unbelievable" estimate~(\ref{pw:f0}) applying for the disk $B_{r_*}(z_0)$ we obtain the required estimate
\begin{equation}
\label{hpf8}
|\we(z_0)-\we_0|\le C'(\delta_k+\e_k). 
\end{equation}

\medskip The remaining two cases ($ii$)--($iii$) are treated similarly using~(\ref{hpf3-}). Because of arbitrariness of $z_0\in \Omega_k$, we proved really that 
\begin{equation}
\label{hpf10}
\sup\limits_{z\in \Omega_k}|\we(z)-\we_0|\to0\qquad\mbox{ as }k\to\infty. 
\end{equation}
By construction, this implies the required uniform convergence~(\ref{eqh1}). Theorem~\ref{half-plane-t} is proved. \qed

\medskip
{\bf \large Proof of Corollary~}\ref{pw-cor1}. Let assumptions of Corollary~\ref{pw-cor1} be fulfilled. Fix $\delta\in(0,1)$. Take arbitrary $z_0\in B_{\delta R}$ and fix it also. 
Denote $r_0=R-|z_0|$. Then by construction
\begin{equation}
\label{dbe-1}
r_0>(1-\delta)R,
\end{equation}
\begin{equation}
\label{dbe-2}
B(z_0,r_0)\subset B_R.
\end{equation}
Denote $\we_0=\bar\we(R)$. Take a ''good'' circle $S_{r_1}$ with radius $r_1\in\biggl(\frac{1+2\delta}3R,\frac{1+\delta}2R\biggr)$ such that 
\begin{equation}
\label{dbe-3}
\sup\limits_{z\in S_{r_1}}|\we(z)-\we_0|\le C_\delta\biggl(\int_{B_R} |\nabla \mathbf{w}|^2\biggr)^\frac12,
\end{equation}
the existence of such ${r_1}$ follows directly from Lemmas~\ref{lem-new4}, \ref{lt3}. By construction, 
$$r_1-|z_0|>\frac{1+2\delta}3R-\delta R=\frac{1-\delta}3R,$$
$$r_0-\bigl(r_1-|z_0|\,\bigr)=R-|z_0|-r_1+|z_0|=R-r_1>\frac{1-\delta}2R,$$
in particular,
\begin{equation}
\label{dbe-4}\frac{r_0}{r_1-|z_0|}>\frac2{1+\delta}.\end{equation}
Now again take a good circle $S_{r_*}(z_0)$ centered at $z_0$ with radius $r_*\in\bigl(r_1-|z_0|,r_0\bigr) $
satisfying 
\begin{equation}
\label{dbe-5}
\sup\limits_{z_1,z_2\in S_{r_*}(z_0)}|\we(z_1)-\we(z_2)|\le C'_\delta\biggl(\int_{B_R} |\nabla \mathbf{w}|^2\biggr)^\frac12,
\end{equation}
the existence of such ${r_*}$ follows directly from (\ref{dbe-4}) and from Lemma~\ref{lt2} applied to the annulus
$$\{z\in\R^2:r_1-|z_0|<|z-z_0|<r_0\}.$$
By construction, since $r_*>r_1-|z_0|$, we have 
$$S_{r_*}(z_0)\cap S_{r_1}\ne\emptyset.$$
Thus, from ''goodness'' of these both circles (see (\ref{dbe-3}) and (\ref{dbe-5})\,) we obtain
\begin{equation}
\label{dbe-6}
\sup\limits_{z\in S_{r_*}(z_0)}|\we(z)-\we_0|\le C''_\delta\biggl(\int_{B_R} |\nabla \mathbf{w}|^2\biggr)^\frac12,
\end{equation}
Then from the "unbelievable" estimate~(\ref{pw:f0}) applying to the ball $B_{r_*}(z_0)$ we obtain the required estimate
$$|\we(z_0)-\we_0|=|\we(z_0)-\bar\we(R)|\le C'''_\delta\biggl(\int_{B_R} |\nabla \mathbf{w}|^2\biggr)^\frac12. 
$$
Because of arbitrariness of $z_0\in B_{\delta R}$, Corollary~\ref{pw-cor1} is proved completely. \qed

\medskip
Corollary~\ref{pw-cor2} can be proved in a similar way, using (\ref{bardiff}) and (\ref{in:estim-m1}) to estimate for the difference between the values~$\bar\we(r_2)$ and $\bar\we(r_1)$. 
\medskip

{\bf Acknowledgement.} {\sl X. Ren is supported by the China Postdoctoral Science Foundation under Grant Number BX20230019 and the National Key R\&D Program of China (No. 2020YFA0712800).}

\

 {\bf Data availability statement.} {\sl Data sharing not applicable to this article as no datasets were generated or analyzed during the current study.}

\

{\bf Conflict of interest statement}. {\sl The authors declare that they have no conflict of interest.}

\bibliographystyle{plain}

\end{document}